\documentclass[11pt]{article} \makeatother
\usepackage{amsfonts,amsmath,amssymb}

\begin{document}

\numberwithin{equation}{section}
\newtheorem{theorem}{\ \ \ \ Theorem}[section]
\newtheorem{proposition}[theorem]{\ \ \ \ Proposition}
\newtheorem{lemma}[theorem]{\ \ \ \ Lemma}
\newtheorem{remark}{\ \ \ \ Remark}
\newcommand{\be}{\begin{equation}}
\newcommand{\ee}{\end{equation}}
\newcommand\bes{\begin{eqnarray}} \newcommand\ees{\end{eqnarray}}
\newcommand{\bess}{\begin{eqnarray*}}
\newcommand{\eess}{\end{eqnarray*}}
\newcommand\ds{\displaystyle}

 \begin{center} {\bf\Large Global nonexistence of solutions for the viscoelastic wave}\\[2mm]
  {\bf\Large   equation of Kirchhoff type with high energy}
\\[4mm]
  {\large Gang Li, Linghui Hong, Wenjun Liu$^{1}$}\\[1mm]
{\small  College of Mathematics and Physics, Nanjing University of
Information Science and Technology, Nanjing 210044, China. E-mail: wjliu@nuist.edu.cn.}\\[1mm]
\end{center}

\setlength{\baselineskip}{17pt}{\setlength\arraycolsep{2pt}

\begin{quote}
\noindent {\bf Abstract:} In this paper we consider the viscoelastic
wave equation of Kirchhoff type:
$$ u_{tt}-M(\|\nabla u\|_{2}^{2})\Delta u+\int_{0}^{t}g(t-s)\Delta u(s){\rm d}s+u_{t}=|u|^{p-1}u $$
 with Dirichlet boundary conditions. Under some suitable
 assumptions on $g$ and the initial data, we established a global nonexistence result for certain solutions with
 arbitrarily high energy.


\end{quote}

\setlength{\baselineskip}{17pt}{\setlength\arraycolsep{2pt}

\section{Introduction}

In this paper we consider the following problem
\bes\left\{\begin{array}{ll} \displaystyle u_{tt}-M(\|\nabla
u\|_{2}^{2})\Delta u+\int_{0}^{t}g(t-s)\Delta u(s){\rm d}
s+|u_{t}|^{m-1}u_{t}=|u|^{p-1}u,\quad &(x,t)\in
\Omega\times(0,\infty),
\medskip\\\medskip
  \displaystyle u(x,t)=0, \quad  &(x,t)\in \partial\Omega\times[0,\infty),\\
\displaystyle u(x,0)=u_{0}(x),\quad u_{t}(x,0)=u_{1}(x), \quad   &x
\in  \Omega, \label{1.1}
 \end{array}\right.
 \ees
where $\Omega$ is a bounded domain in $\mathbb{R}^n$ ($n\geq 1$)
with a  smooth boundary $\partial\Omega,$ $p>1,$ $M(s)$ is a
nonnegative $C^{1}$ function like $M(s)=a+bs^{\gamma}$ for $s\geq0$,
$a\geq0, b\geq0, a+b>0, \gamma>0$ and $g(t)$
 represents the kernel of memory term.

  Problem \eqref{1.1} without the viscoelastic term ($i.e., g=0$)
has been extensively studied and many results concerning global
existence, decay and blow-up have been established. For example, the
following equation
\begin{equation}
u_{tt}-M(\|\nabla u\|_{2}^{2})\Delta u+g(u_{t})=f(u), \quad (x,t)\in
\Omega\times(0,\infty), \label{1.2}
\end{equation}
has been considered by Matsuyama and Ikehata in \cite{33b}, for
$g(u_{t})=\delta|u_{t}|^{p-1}u_{t}$ and $f(u)=\mu|u|^{q-1}u.$ The
authors proved existence of the global solutions by using
Faedo-Galerkin method and the decay of energy based on the method of
Nakao \cite{2005a}--\cite{203a}. Later, Ono \cite{3d} investigated
equation \eqref{1.2} for $M(s)=bs^{\gamma}$ and $f(u)=|u|^{p-2}u.$
When $g(u_{t})=-\Delta
 u_{t},$ $u_{t}$ or $|u_{t}|^{\beta}u_{t},$ the author showed that the solutions blow up in finite
 time with $E(0)\leq 0.$ For $M(s)=a+bs^{\gamma}$ and $g(u_{t})=u_{t},$ this model was considered by
 the same author in \cite{8d}. By applying the potential well
 method he obtained the blow-up properties with positive initial energy $E(0).$ Recently,
  Zeng et al. \cite {22b} studied equation \eqref{1.2} for the case $g(u_{t})=u_{t}$ with initial condition and zero Dirichlet boundary
condition.
 By using the concavity argument, they proved
 that the solutions to
equation \eqref {1.2} blow up in finite time with arbitrarily high
energy.

 In the case of $M\equiv1$ and in the presence of the viscoelastic term ($i.e., g\neq0$), the equation
 \begin{equation}
u_{tt}-\Delta u+\int_{0}^{t}g(t-s)\Delta u(s){\rm
d}s+|u_{t}|^{m-1}u_{t}=|u|^{p-1}u, \quad (x,t)\in \Omega\times(0,
\infty), \label{1.3}
  \end{equation}
  was studied by Messaoudi in \cite{18e}, where the author proved that any weak
  solution with negative initial energy blows up in finite time if
  $p>m$ and $$\int_{0}^{\infty}g(s){\rm d}s\leq\frac{p-1}{p-1+1/(p+1)},$$
  while the solution continue to exist globally for any initial data in the appropriate
  space if $m\geq p.$ This blow-up result was improved by the same
  author in \cite{19e} for positive initial energy under suitable
  conditions on $g,$ $m$ and $p.$ More recently, Wang \cite{11a}
investigated equation \eqref {1.3} and established a blow-up result
   with arbitrary positive initial energy. In the related work, Cavalcanti et al.
 \cite {5e} studied the following equation
 \begin{equation}
u_{tt}-\Delta u+\int_{0}^{t}g(t-s)\Delta u(s){\rm
d}s+a(x)u_{t}+|u|^{\gamma}u=0,\quad (x,t)\in \Omega\times(0,
\infty), \label{1.4}
 \end{equation}
 where $a:\Omega\longrightarrow R^{+}$ is a function which may be
 null on a part of $\Omega.$ Under the condition that $a(x)\geq
 a_{0}>0$ on $\omega\subset\Omega,$ with $\omega$ satisfying some
 geometric restrictions and $-\xi_{1}g(t)\leq g'(t)\leq
 -\xi_{2}g(t),$ $
 t\geq0$ to guarantee $\|g\|_{L^{1}((0, \infty))}$ is small enough,
 they proved an exponential decay rate.

 When $g\neq0$ and $M$ is not a constant function, problems related
 to \eqref {1.1} have been treated by many authors. Wu
 and Tsai \cite{34c} considered the global existence, asymptotic behavior and blow-up properties for the following
 equation
\begin{equation}
u_{tt}-M(\|\nabla u\|_{2}^{2})\Delta u+\int_{0}^{t}g(t-s)\Delta
u(s){\rm d}s-\Delta u_{t}=f(u), \quad (x,t)\in \Omega\times(0,
\infty), \label{1.5}
\end{equation}
with the same initial and boundary conditions as that of problem
\eqref{1.1}. To obtain the decay result, they assumed that the
nonnegative kernel $g'(t)\leq -rg(t), \forall\ t\geq0$ for some
$r>0.$ In \cite{33c}, Wu then extended the result of \cite{34c}
under a weaker condition on $g$ (i.e., $ g'(t)\leq 0$ for $t\geq0$).
For other papers related to existence, uniform decay and blow-up of
solutions of nonlinear wave equations, see \cite{c2007, c2004,
2001a, 2002a, 2003a, 13a, wljy, 9d, 10d, sw2006, 2004a, 32b, y2009}
and references therein.

Motivated by the above research, we consider problem \eqref{1.1} for
$m=1$ in this paper and  establish a global nonexistence result for
certain solutions with arbitrarily high energy. In this way, we can
extend the result of \cite{22b} to nonzero term $g$ and the result
of \cite{11a} to nonconstant $M(s)$. We also obtain the new result
for blow-up properties of local solution with arbitrarily high
energy. Throughout the rest of this paper, we always assume that
$m=1.$

 The structure of this
paper is as follows. In section 2, we present some assumptions,
notations and main result. Section 3 is devoted to the proof of the
main result.
\section{Preliminaries and main result}

In this section, we shall give some assumptions, notations and main
result. We first give the following assumptions:
\begin{quote}
(A1) $g\in C^{1}([0, \infty))$ is a non-negative and non-increasing
function satisfying
$$1-\int_{0}^{\infty}g(s){\rm d}s=l>0.$$
(A2) The function $e^{\frac{t}{2}}g(t)$ is of positive type in the
following sense:
$$\int_{0}^{t}v(s)\int_{0}^{s}e^{\frac{s-z}{2}}g(s-z)v(z){\rm d}z{\rm d}s\geq0,
\quad \forall\ v\in C^{1}([0, \infty)) \quad \text{and}\quad
\forall\ t>0.$$ In order to prove our result, we make the following assumption on $M$ and $g:$ \\
 (A3) There exist two positive constants,
$m_{1}$ and $\alpha,$ such that
$$\frac{p+1}{2}\overline{M}(s)-\left[M(s)+\frac{p+1}{2}\int_{0}^{t}g(\tau){\rm d}\tau\right]s\geq m_{1}s^{\alpha},
 \quad \forall\ s\geq0,$$ where
 $\displaystyle\overline{M}(s)=\int_{0}^{s}M(\tau){\rm d}\tau.$
\end{quote}

\begin{remark} \label{re1}
 It is clear that when $M(s)=a+bs^{\gamma}$
 for $s\geq0$, $a\geq0,$ $b\geq0,$ $a+b>0,$ $\gamma>0$ and $p>1+2\gamma$, condition
 $($A3$)$ can be replaced by
 \begin{align}
 \int_{0}^{\infty}g(\tau){\rm d}\tau<\left\{ {\begin{array}{l}
 \frac{p-1}{p+1}a,\quad \quad \quad \quad \quad \quad \quad \quad \text{if}\quad
  a>0 \quad \text{and}\quad b\geq0, \\
 \frac{(p-1-2\gamma)b}{C_{p}^{\gamma}(p+1)(\gamma+1)}\|u_{0}\|_{2}^{2\gamma}, \quad \quad\text{if}\quad
 a=0 \quad \text{and} \quad b>0, \\
\end{array}}\right. \label{*}
\end{align}
 which is the same as the one in \cite [Theorem
1.1]{11a} for the case $a=1$ and
$b=0$, where $C_{p}$ is the constant from the Poincar\'e inequality $\|u(t)\|_{2}^2\le C_p \|\nabla u(t)\|_{2}^2$.

Indeed, by straightforward calculation, we obtain
\begin{align*}
&\frac{p+1}{2}\overline{M}(s)-\left[M(s)+\frac{p+1}{2}\int_{0}^{t}g(\tau){\rm
d}\tau\right]s
\\
=&\frac{p+1}{2}\left(as+\frac{b}{\gamma+1}s^{\gamma+1}\right)-as-
bs^{\gamma+1}-\frac{(p+1)s}{2}\int_{0}^{t}g(\tau){\rm d}\tau \\
=&\frac{p-1}{2}as+\frac{(p-1-2\gamma)b}{2(\gamma+1)}s^{\gamma+1}-\frac{(p+1)s}{2}\int_{0}^{t}g(\tau){\rm
d}\tau.
\end{align*}

If $a>0$ and $b\geq0,$ it follows from  \eqref{*} that
$\int_{0}^{\infty}g(\tau){\rm d}\tau<\frac{p-1}{p+1}a.$  Thus, we have
\begin{align*}
&\frac{p+1}{2}\overline{M}(s)-\left[M(s)+\frac{p+1}{2}\int_{0}^{t}g(\tau){\rm
d}\tau\right]s
\\
>&\frac{p-1}{2}as+\frac{(p-1-2\gamma)b}{2(\gamma+1)}s^{\gamma+1}-\frac{(p+1)s}{2}\left[\frac{p-1}{p+1}a
-\frac{\frac{p-1}{p+1}a-\int_{0}^{\infty}g(\tau){\rm
d}\tau}{2}\right]\\
=&\frac{(p-1-2\gamma)b}{2(\gamma+1)}s^{\gamma+1}+\frac{1}{2}\left[\frac{p-1}{p+1}a-\int_{0}^{\infty}g(\tau){\rm
d}\tau\right]s
\geq\frac{1}{2}\left[\frac{p-1}{p+1}a-\int_{0}^{\infty}g(\tau){\rm
d}\tau\right]s.
\end{align*}
Therefore, we can choose
$m_{1}=\frac{1}{2}\left[\frac{p-1}{p+1}a-\int_{0}^{\infty}g(\tau){\rm
d}\tau\right]$ and $\alpha=1$ in  condition
 $($A3$)$. \newpage

If $a=0$ and $b>0,$ then
\begin{align*}
&\frac{p+1}{2}\overline{M}(s)-\left[M(s)+\frac{p+1}{2}\int_{0}^{t}g(\tau){\rm
d}\tau\right]s=\frac{(p-1-2\gamma)b}{2(\gamma+1)}s^{\gamma+1}-\frac{(p+1)s}{2}\int_{0}^{t}g(\tau){\rm
d}\tau \\
>&\frac{(p-1-2\gamma)b}{2(\gamma+1)}s^{\gamma+1}-\frac{(p+1)s}{2}\left[\frac{(p-1-2\gamma)b}{C_{p}^{\gamma}(p+1)(\gamma+1)}
\|u_{0}\|_{2}^{2\gamma}-\frac{\frac{(p-1-2\gamma)b}{C_{p}^{\gamma}(p+1)(\gamma+1)}\|u_{0}\|_{2}^{2\gamma}-
\int_{0}^{\infty}g(\tau){\rm d}\tau}{2}\right]\\
=&\frac{(p-1-2\gamma)b}{2(\gamma+1)}s\left(s^{\gamma}-\frac{1}{C_{p}^{\gamma}}\|u_{0}\|_{2}^{2\gamma}\right)
+\frac{(p+1)s}{4}\left[\frac{(p-1-2\gamma)b}{C_{p}^{\gamma}(p+1)(\gamma+1)}\|u_{0}\|_{2}^{2\gamma}-
\int_{0}^{\infty}g(\tau){\rm d}\tau \right].
\end{align*}
Taking $s=\|\nabla u(t)\|_{2}^{2},$ applying Lemma \ref{3.3} below and
Poincar\'e's inequality, we can get
\begin{align*}
&\frac{p+1}{2}\overline{M}(\|\nabla u(t)\|_{2}^{2})-\left[M(\|\nabla
u(t)\|_{2}^{2})+\frac{p+1}{2}\int_{0}^{t}g(\tau){\rm
d}\tau\right]\|\nabla
u(t)\|_{2}^{2} \\
>&\frac{(p-1-2\gamma)b}{2(\gamma+1)}\|\nabla u(t)\|_{2}^{2}\left(\|\nabla u(t)
\|_{2}^{2\gamma}-\frac{1}{C_{p}^{\gamma}}\|u_{0}\|_{2}^{2\gamma}\right)\\
&+\frac{p+1}{4}\|\nabla
u(t)\|_{2}^{2}\left[\frac{(p-1-2\gamma)b}{C_{p}^{\gamma}(p+1)(\gamma+1)}\|u_{0}\|_{2}^{2\gamma}-
\int_{0}^{\infty}g(\tau){\rm d}\tau \right]\\
\geq &\frac{(p-1-2\gamma)b}{2(\gamma+1)}\|\nabla
u(t)\|_{2}^{2}\left(\|\nabla u(t)
\|_{2}^{2\gamma}-\frac{1}{C_{p}^{\gamma}}\|u(t)\|_{2}^{2\gamma}\right)\\
&+\frac{p+1}{4}\|\nabla
u(t)\|_{2}^{2}\left[\frac{(p-1-2\gamma)b}{C_{p}^{\gamma}(p+1)(\gamma+1)}\|u_{0}\|_{2}^{2\gamma}-
\int_{0}^{\infty}g(\tau){\rm d}\tau \right]\\
\geq&\frac{(p-1-2\gamma)b}{2(\gamma+1)}\|\nabla
u(t)\|_{2}^{2}\left(\|\nabla u(t)
\|_{2}^{2\gamma}-\|\nabla u(t)\|_{2}^{2\gamma}\right)\\
&+\frac{p+1}{4}\|\nabla
u(t)\|_{2}^{2}\left[\frac{(p-1-2\gamma)b}{C_{p}^{\gamma}(p+1)(\gamma+1)}\|u_{0}\|_{2}^{2\gamma}-
\int_{0}^{\infty}g(\tau){\rm d}\tau \right]\\
=&\frac{p+1}{4}\left[\frac{(p-1-2\gamma)b}{C_{p}^{\gamma}(p+1)(\gamma+1)}\|u_{0}\|_{2}^{2\gamma}-
\int_{0}^{\infty}g(\tau){\rm d}\tau \right]\|\nabla u(t)\|_{2}^{2}.
\end{align*}
So, we can choose
$m_{1}=\frac{p+1}{4}\left[\frac{(p-1-2\gamma)b}{C_{p}^{\gamma}(p+1)(\gamma+1)}\|u_{0}\|_{2}^{2\gamma}-
\int_{0}^{\infty}g(\tau){\rm d}\tau \right]$ and $\alpha=1$ in  condition
 $($A3$)$.
\end{remark}

Next, we introduce some notations. The energy functional $E(t)$ and
an auxiliary functional $I(u)$ of the solution $u(t)$ of problem
\eqref {1.1} are defined as follows:
\begin{equation}
E(t):=E(u(t))=\frac{1}{2}\|u_{t}\|_{2}^{2}+\frac{1}{2}\overline{M}(\|\nabla
u\|_{2}^{2})-\frac{1}{2}\int_{0}^{t}g(s){\rm d}s\|\nabla
u\|_{2}^{2}+\frac{1}{2}(g\circ \nabla
u)(t)-\frac{1}{p+1}\|u\|_{p+1}^{p+1},\label{2.1}
\end{equation}
 and
\begin{equation}
 I(u)=M(\|\nabla u\|_{2}^{2})\|\nabla
u\|_{2}^{2}-\|u\|_{p+1}^{p+1},\label {2.2}
\end{equation}
where $$(g\circ
w)(t)=\int_{0}^{t}g(t-s)\|w(t,\cdot)-w(s,\cdot)\|_{2}^{2}{\rm d}s.$$

As in \cite{11a, 22b}, we can get
\begin{equation}
\frac{d}{{\rm d}t}E(t)=-\|u_{t}\|_{2}^{2}-\frac{1}{2}g(t)\|\nabla
u\|_{2}^{2}+\frac{1}{2}(g'\circ \nabla u)(t)\leq0,\label {2.3}
\end{equation}
for $t\geq0$. Then we have
\begin{equation}
E(t)=E(0)-\int_{0}^{t}\|u_{s}\|_{2}^{2}{\rm
d}s+\frac{1}{2}\int_{0}^{t}(g'\circ \nabla
u)(s){\rm d}s\\
-\frac{1}{2}\int_{0}^{t}g(s) \|\nabla u(s)\|_{2}^{2}{\rm d}s.\label
{2.4}
\end{equation}

Now we are in a position to state our main result.

\begin{theorem}\label{Th2.1}
Assume that $M$ and $g$ satisfy assumptions $($A1$)$-$($A3$)$.
Suppose further that $1<p\leq\frac{n}{n-2}$ when $n\geq3,$
$1<p<\infty$ when $n=1,2.$ Let $u$ be a solution of problem
\eqref{1.1} with initial data $u_{0}\in H_{0}^{1}(\Omega)\cap
H^{2}(\Omega),$ $u_{1}\in H_{0}^{1}(\Omega)$ satisfying
\begin{align}
E(0)&>0,\label {2.5}
\\
I(u_{0})&<0,\label {2.6}
\\
\int_{\Omega}u_{0}u_{1}{\rm d}x&>0,\label {2.7}
\\
\|u_{0}\|_{2}^{2}&>C_{p}\left(\frac{p+1}{m_{1}}E(0)\right)^{1/\alpha}.\label
{2.8}
\end{align}
 Then the solution of problem \eqref{1.1} blows up
in finite time $0<T^{*}<+\infty,$ which means that
\begin{align}
\lim_{t\rightarrow
T^{*-}}\left(\|u(t)\|_{2}^{2}+\int_{0}^{t}\|u(s)\|_{2}^{2}{\rm
d}s\right)=\infty,\label {2.9}
\end{align}
where $C_{p}$ is a constant from the Poincar\'e inequality and
$m_{1}$ comes from condition
 $($A3$)$.
\end{theorem}

\section{Proof of main result}

Before we start to prove Theorem \ref{Th2.1}, it is necessary to
state the local existence theorem for problem \eqref{1.1}, whose
proof follows the arguments in \cite{10d,34c}.

\begin{theorem}\label{Th3.1}
Assume that $($A1$)$ holds, and $1<p\leq\frac{n}{n-2}$ when
$n\geq3,$ $1<p<\infty$ when $n=1,2.$ For $u_{0}\in
H_{0}^{1}(\Omega)\cap H^{2}(\Omega)$, $u_{1}\in H_{0}^{1}(\Omega)$
and $M(\|\nabla u_{0}\|_{2}^{2})>0$, problem \eqref{1.1} has a
unique local solution
$$u\in C([0,T]; H_{0}^{1}(\Omega)\cap
H^{2}(\Omega)),\quad u_{t}\in C([0,T]; L^{2}(\Omega))\cap
L^{2}([0,T]; H_{0}^{1}(\Omega))$$ for the maximum existence time
$T>0.$
\end{theorem}

The proof of Theorem \ref{Th2.1} relies on the following lemmas.

\begin{lemma}{\rm(\cite [Lemma 2.1]{11a})} \label{le3.2}
Assume that $g(t)$ satisfies assumptions $($A1$)$-$($A2$)$, and
$H(t)$ is a function which is twice continuously differentiable
satisfying \bes\left\{\begin{array}{ll}
H''(t)+H'(t)>\int_{0}^{t}g(t-s)\int_{\Omega}\nabla u(s)\nabla u(t)d
x d s\medskip\\\medskip
  \displaystyle H(0)>0,\quad H'(0)>0,  \label{3.1}
\end{array}\right.
 \ees
 for every $t\in[0,T),$ where $u(t)$ is the corresponding
 solution of problem \eqref{1.1} with $u_{0}$ and $u_{1}.$ Then the
 function $H(t)$ is strictly increasing on $[0,T).$
\end{lemma}

\begin{lemma} \label{le3.3}
Suppose that $u_{0}\in H_{0}^{1}(\Omega)\cap H^{2}(\Omega)$ and
$u_{1}\in H_{0}^{1}(\Omega)$ satisfy
\begin{align}
\int_{\Omega}u_{0}u_{1}{\rm d}x>0.\label{3.2}
\end{align}
If the solution $u(t)$ of problem \eqref{1.1} exists on $[0,T)$ and
satisfies
\begin{align}
I(u(t))<0,\label{3.3}
\end{align}
then $\|u(t)\|_{2}^{2}$ is strictly increasing on $[0,T).$
\end{lemma}

{\bf Proof.}\  Since $u(t)$ is the solution of problem \eqref{1.1},
by a simple computation, we have
\begin{align*}
\frac{1}{2}\frac{d^{2}}{{\rm d}t^{2}}\int_{\Omega}|u(x,t)|^{2}{\rm d}x=&\int_{\Omega}(|u_{t}|^{2}+uu_{tt}){\rm d}x\\
=&\|u_{t}\|_{2}^{2}-M(\|\nabla u\|_{2}^{2})\|\nabla
u\|_{2}^{2}+\|u\|_{p+1}^{p+1}\\
 &+\int_{0}^{t}g(t-s)\int_{\Omega}\nabla u(s)\nabla
 u(t){\rm d}x{\rm d}s-\int_{\Omega}uu_{t}{\rm d}x\\
 >&-\int_{\Omega}uu_{t}{\rm d}x+\int_{0}^{t}g(t-s)\int_{\Omega}\nabla u(s)\nabla
 u(t){\rm d}x{\rm d}s,
\end{align*}
where the last inequality is derived by \eqref{3.3}.
Then we get
\begin{align*}
\frac{d^{2}}{{\rm d}t^{2}}\int_{\Omega}|u(x,t)|^{2}{\rm
d}x+\frac{d}{{\rm d}t}\int_{\Omega}|u(x,t)|^{2}{\rm
d}x>\int_{0}^{t}g(t-s)\int_{\Omega}\nabla u(s)\nabla
 u(t){\rm d}x{\rm d}s.
\end{align*}
Therefore, by using Lemma \ref{le3.2}, we finish our proof. $\Box$

\begin{lemma}  \label{le3.4}
If $u_{0}\in H_{0}^{1}(\Omega)\cap H^{2}(\Omega)$ and $u_{1}\in
H_{0}^{1}(\Omega)$ satisfy the assumptions in Theorem \ref{Th2.1},
then the solution $u(t)$ of problem \eqref{1.1} satisfies
\begin{align}
I(u(t))&<0,\label{3.4}
\\
\|u\|_{2}^{2}&>C_{p}\left(\frac{p+1}{m_{1}}E(0)\right)^{1/\alpha},\label{3.5}
\end{align}
for all $t\in [0,T).$
\end{lemma}

{\bf Proof.}\ We will prove the above lemma by contradiction. First
we assume that \eqref{3.4} is not true over $[0,T),$ it means that
there exists a time $t_{0}$ such that
\begin{align}
t_{0}=\min\{t\in (0,T): I(u(t))=0 \}.\label{3.6}
\end{align}
Since $I(u(t))<0$ on $[0,t_{0}),$ by Lemma \ref{le3.3}, we see that
$\displaystyle\int_{\Omega}u^{2}{\rm d}x$ is strictly increasing
over $[0,t_{0}),$ which implies
\begin{align}
\int_{\Omega}u^{2}{\rm d}x>\int_{\Omega}u_{0}^{2}{\rm
d}x>C_{p}\left(\frac{p+1}{m_{1}}E(0)\right)^{1/\alpha}.\label{3.7}
\end{align}
And by the continuity of $\displaystyle\int_{\Omega}u^{2}{\rm d}x$
on $t,$ we note that
\begin{align}
\int_{\Omega}u^{2}(t_{0}){\rm
d}x>C_{p}\left(\frac{p+1}{m_{1}}E(0)\right)^{1/\alpha}.\label{3.8}
\end{align}
On the other hand, by \eqref{2.1} and \eqref{2.4}, we get
\begin{align}
\overline{M}(\|\nabla u(t_{0})\|_{2}^{2})-\int_{0}^{t_{0}}g(s){\rm
d}s\|\nabla u(t_{0})\|_{2}^{2}+(g\circ \nabla
u)(t_{0})-\frac{2}{p+1}\|u(t_{0})\|_{p+1}^{p+1}
 \leq2E(0).\label{3.9}
\end{align}
Combining \eqref{3.9} with \eqref{3.6} yields
\begin{align}
 \frac{p+1}{2}\overline{M}&(\|\nabla
u(t_{0})\|_{2}^{2})-\frac{p+1}{2}\int_{0}^{t_{0}}g(s){\rm
d}s\|\nabla
u(t_{0})\|_{2}^{2} \nonumber \\
 &+\frac{p+1}{2}(g\circ \nabla
u)(t_{0}) -M(\|\nabla u(t_{0})\|_{2}^{2})\|\nabla
u(t_{0})\|_{2}^{2}\leq(p+1)E(0).
\label{3.10}
\end{align}
By (A3), we get
\begin{align}
m_{1}\|\nabla u(t_{0})\|_{2}^{2\alpha}<(p+1)E(0). \label{3.11}
\end{align}
i.e.,
\begin{align}
\|\nabla
u(t_{0})\|_{2}^{2}<\left(\frac{p+1}{m_{1}}E(0)\right)^{1/\alpha}.
\label{3.12}
\end{align}
By Poincar\'e's inequality, we have
\begin{align}
\|u(t_{0})\|_{2}^{2}<C_{p}\left(\frac{p+1}{m_{1}}E(0)\right)^{1/\alpha}.
\label{3.13}
\end{align}

Obviously, there is a contradiction between  \eqref{3.8} and
\eqref{3.13}, thus we prove that
\begin{align}
I(u(t))<0,\label{3.14}
\end{align}
for every $t\in(0,T).$ By Lemma \ref{le3.3}, it follows that
$\displaystyle\int_{\Omega}u^{2}{\rm d}x$ is strictly increasing on
$[0,T),$ which implies that
\begin{align}
\int_{\Omega}u^{2}{\rm d}x\geq\int_{\Omega}u_{0}^{2}{\rm
d}x>C_{p}\left(\frac{p+1}{m_{1}}E(0)\right)^{1/\alpha},\label{3.15}
\end{align}
for every $t\in[0,T).$ This completes the proof of Lemma \ref{le3.4}
. $\Box$

\begin{lemma}\rm{(\cite{2006c})}  \label{le3.5}
Assume that  $P(t)\in C^{2},$ $P(t)\geq0,$ satisfies the inequality
\begin{align*}
P(t)P''(t)-(1+\theta)P'^{2}(t)\geq0,
\end{align*}
for certain real number $\theta>0,$ and $P(0)>0,$ $P'(0)>0.$ Then
there exists a real number $T^{*}$ such that $0<T^{*}\leq
P(0)/\theta P'(0)$ and

$$P(t)\rightarrow \infty$$ as $t\rightarrow
T^{*-}.$

\end{lemma}

{\bf Proof of Theorem \ref{Th2.1}.}\ We prove our main result by
adopting concavity method, and define an auxiliary function by
\begin{align}
G(t)=\|u(t)\|_{2}^{2}+\int_{0}^{t}\|u(s)\|_{2}^{2}{\rm
d}s+(T_{0}-t)\|u_{0}\|_{2}^{2}+\beta(t_{2}+t)^{2},\label{3.16}
\end{align}
where $T_{0}, t_{2}$ and $\beta$ are positive constants, which will
be chosen later.

A straightforward calculation gives
\begin{align}
G'(t)&=2\int_{\Omega}uu_{t}{\rm
d}x+\|u(t)\|_{2}^{2}-\|u_{0}\|_{2}^{2}+2\beta(t_{2}+t)
\nonumber\\
&=2\int_{\Omega}uu_{t}{\rm d}x+2\int_{0}^{t}(u(s),u_{s}(s)){\rm
d}s+2\beta(t_{2}+t), \label{3.17}
\end{align}
consequently,
\begin{align}
G''(t)=&2\int_{\Omega}|u_{t}|^{2}{\rm d}x+2\int_{\Omega}uu_{tt}{\rm
d}x+2\int_{\Omega}uu_{t}{\rm d}x+2\beta
\nonumber \\
 =&2\|u_{t}\|_{2}^{2}-2M(\|\nabla u\|_{2}^{2})\|\nabla
u\|_{2}^{2}+2\int_{0}^{t}g(t-s)\int_{\Omega}\nabla u(s)\nabla
u(t){\rm d}x{\rm d}s \nonumber \\
&-2\int_{\Omega}uu_{t}{\rm d}x+2\int_{\Omega}uu_{t}{\rm
d}x+2\|u\|_{p+1}^{p+1}+2\beta
\nonumber \\
 =&2\|u_{t}\|_{2}^{2}-2M(\|\nabla u\|_{2}^{2})\|\nabla
u\|_{2}^{2}+2\|u\|_{p+1}^{p+1}+2\int_{0}^{t}g(t-s){\rm d}s\|\nabla
u\|_{2}^{2} \nonumber \\&+2\int_{0}^{t}g(t-s)\int_{\Omega}\nabla
u(t)(\nabla u(s)-\nabla u(t)){\rm d}x{\rm d}s+2\beta. \label{3.18}
\end{align}
We will use the following Young inequality to estimate the fifth
term of the right hand side of \eqref{3.18},
\begin{align*}
rs\leq\frac{r^{2}}{2\epsilon}+\frac{\epsilon s^{2}}{2}
\end{align*}
where $\epsilon=\frac{1}{2},$ $r\geq0$ and $s\geq0.$ We obtain
\begin{align}
 \int_{0}^{t}g(t-s)\int_{\Omega}|\nabla u(t)||\nabla u(s)-\nabla
u(t)|{\rm d}x{\rm d}s\leq\int_{0}^{t}g(s){\rm d}s\|\nabla
u(t)\|_{2}^{2}+\frac{1}{4}(g\circ \nabla u)(t), \label{3.19}
\end{align}

Substitute \eqref{2.1} and \eqref{3.19} for the third and the fifth
terms of the right hand side of \eqref{3.18}, respectively, we have
\begin{align}
G''(t)\geq & (p+3)\|u_{t}\|_{2}^{2}+(p+1)\overline{M}(\|\nabla
u\|_{2}^{2})-2M(\|\nabla u\|_{2}^{2})\|\nabla
u\|_{2}^{2}-(p+1)\int_{0}^{t}g(s){\rm d}s\|\nabla u\|_{2}^{2} \nonumber \\
&-2(p+1)E(t)+\left(p+\frac{1}{2}\right)(g\circ \nabla
u)(t)+2\beta.\label{3.20}
\end{align}
By (A3), we deduce
\begin{align}
G''(t)>(p+3)\|u_{t}\|_{2}^{2}+2m_{1}\|\nabla
u\|_{2}^{2\alpha}-2(p+1)E(t)+\left(p+\frac{1}{2}\right)(g\circ
\nabla u)(t)+2\beta.\label{3.21}
\end{align}
Noting that  \eqref{2.4}, we obtain that
\begin{align}
 -E(t)\geq-E(0)+\int_{0}^{t}\|u_{s}\|_{2}^{2}{\rm d}s. \label{3.22}
\end{align}
Combining \eqref{3.21}-\eqref{3.22} yields
\begin{align}
G''(t)>&(p+3)\|u_{t}\|_{2}^{2}+2m_{1}\|\nabla u\|_{2}^{2\alpha}
-2(p+1)E(0)+\left(p+\frac{1}{2}\right)(g\circ \nabla u)(t) \nonumber \\
&+2(p+1)\int_{0}^{t}\|u_{s}\|_{2}^{2}{\rm d}s+2\beta, \label{3.23}
\end{align}
by \eqref{3.5}, we see that
\begin{align*}
2m_{1}\|\nabla
u\|_{2}^{2\alpha}-2(p+1)E(0)+\left(p+\frac{1}{2}\right)(g\circ
\nabla u)(t)>0.
\end{align*}
 which means that $G''(t)>0$ for every
$t\in(0,T).$
 Thus, by $G'(0)>0$ and $G(0)>0,$
we get $G'(t)$ and $G(t)$ are strictly increasing on $[0,T).$

Thus, we can let $\beta$ satisfy
\begin{align*}
(p+1)\beta<2m_{1}\|\nabla
u\|_{2}^{2\alpha}-2(p+1)E(0)+(p+\frac{1}{2})(g\circ \nabla u)(t).
\end{align*}
Moreover, we let $T_{0}$ and $t_{2}$ satisfy that
$$T_{0}\geq\frac{4}{p-1}\frac{G(0)}{G'(0)},$$
$$\frac{p-1}{2}\left(\int_{\Omega}u_{0}u_{1}{\rm d}x+\beta t_{2}\right)\geq\|u_{0}\|_{2}^{2}.$$
Letting
\begin{align*}
A:=&\|u(t)\|_{2}^{2}+\int_{0}^{t}\|u(s)\|_{2}^{2}{\rm
d}s+\beta(t_{2}+t)^{2},
\\
B:=&\frac{1}{2}G'(t), \\
C:=&\|u_{t}(t)\|_{2}^{2}+\int_{0}^{t}\|u_{s}(s)\|_{2}^{2}{\rm
d}s+\beta.
\end{align*}

Since we have assumed that the solution $u(t)$ to problem
\eqref{1.1} exists for every $t\in[0,T),$ where $T$ is sufficiently
large, we have
\begin{align*}
G(t)&\geq A, \\
G''(t)&\geq (p+3)C
\end{align*}
for every $t\in[0,T_{0}).$ Then it follows that
\begin{align*}
 G''(t)G(t)-\frac{p+3}{4}(G'(t))^{2}\geq(p+3)(AC-B^{2}).
\end{align*}
Furthermore, we have
\begin{align*}
Ar^{2}-2Br+C=\int_{\Omega}(ru(t)-u_{t}(t))^{2}{\rm d}x
+\int_{0}^{t}\|ru(s)-u_{s}(s)\|_{2}^{2}{\rm
d}s+\beta[r(t_{2}+t)-1]^{2}\geq0,
\end{align*}
for every $r\in\mathbb{R},$ which implies that $B^{2}-AC\leq0.$

Thus, we obtain
\begin{align*}
G''(t)G(t)-\frac{p+3}{4}(G'(t))^{2}\geq0,
\end{align*}
for every $t\in[0,T).$

As $\frac{p+3}{4}>1$, letting $\theta=\frac{p-1}{4},$ according to
Lemma \ref{le3.5}, there exists a real number $T^{*}$ such that
$T^{*}<G(0)/\theta G'(0)\leq T_{0}$ and we have
\begin{align*}
\lim_{t\rightarrow T^{*-}}G(t)=\infty,
\end{align*}
i.e.,
 \begin{align}
\lim_{t\rightarrow
T^{*-}}\left(\|u(t)\|_{2}^{2}+\int_{0}^{t}\|u(s)\|_{2}^{2}{\rm
d}s\right)=\infty.\label{3.24}
\end{align}
This completes the proof of Theorem \ref{Th2.1}. $\Box$

\subsection*{Acknowledgments}
This work was partly supported by the Tianyuan Fund of Mathematics (Grant No. 11026211)
and the Natural Science Foundation of the Jiangsu Higher Education Institutions (Grant No.
09KJB110005).

\end{document}